\documentclass{article}

\usepackage{PRIMEarxiv}

\usepackage[utf8]{inputenc} 
\usepackage[T1]{fontenc}    
\usepackage{hyperref}       
\usepackage{url}            
\usepackage{booktabs}       
\usepackage{amsfonts}       
\usepackage{nicefrac}       
\usepackage{microtype}      
\usepackage{lipsum}
\usepackage{fancyhdr}       
\usepackage{graphicx}       
\graphicspath{{media/}}     

\title{Stabilization of stochastic dynamical systems of a random structure with Markov switches and Poisson perturbations
\thanks{\textit{\underline{Citation}}: 
\textbf{Lukashiv, T.; Litvinchuk, Y.; Malyk, I.V.; Golebiewska, A.; Nazarov,
P.V. Stabilization of stochastic dynamical systems of a random structure with Markov switches and Poisson perturbations. Pages.... DOI:000000/11111.}} 
}

\author{
  Taras Lukashiv \\
  Multiomics Data Science Research Group, NORLUX Neuro-Oncology Laboratory \\
  Department of Cancer Research \\
  Luxembourg Institute of Health, L-1210, Luxembourg, Luxembourg\\
  Department of  Mathematical Problems of Control and Cybernetics,\\ 
  Yuriy Fedkovych Chernivtsi National University, 58000, Chernivtsi, Ukraine\\
  \texttt{t.lukashiv@gmail.com} \\
   \And
  Yuliia Litvinchuk,  Igor V. Malyk\\
  Department of  Mathematical Problems of Control and Cybernetics,\\ 
  Yuriy Fedkovych Chernivtsi National University, 58000, Chernivtsi, Ukraine\\
  \texttt{y.litvinchuk@chnu.edu.ua, i.malyk@chnu.edu.ua} \\
  \And
  Anna Golebiewska \\
  NORLUX Neuro-Oncology Laboratory,  Department of Cancer Research \\
  Luxembourg Institute of Health, L-1210, Luxembourg, Luxembourg\\
  \texttt{anna.golebiewska@lih.lu} \\
  \And
  Petr. V. Nazarov \\
  Multiomics Data Science Research Group,  Department of Cancer Research \\
  Luxembourg Institute of Health, L-1445, Strassen, Luxembourg\\
  \texttt{petr.nazarov@lih.lu} \\
}

\begin{document}
\maketitle

\begin{abstract}
An optimal control for a dynamical system optimizes a certain objective function. Here we consider the construction of an optimal control for a stochastic dynamical system with a random structure, Poisson perturbations and random jumps, which makes the system stable in probability. Sufficient conditions of the stability in probability are obtained, using the second Lyapunov method, in which the construction of the corresponding functions plays an important role. Here we provide a solution to the problem of optimal stabilization in a general case. For a linear system with a quadratic quality function, we give a method of synthesis of optimal control based on the solution of Riccati equations. Finally, in an autonomous case, a system of differential equations was constructed to obtain unknown matrices that are used for the building of an optimal control. The method of a small parameter is justified for the algorithmic search of an optimal control. This approach brings a novel solution to the problem of optimal stabilization for a stochastic dynamical system with a random structure, Markov switches and Poisson perturbations.
\end{abstract}


\section{Introduction}
The main problem considered in this paper is the problem of synthesis of optimal control for a controlled dynamical system, described by a stochastic differential equation (SDE) with Poisson perturbations and external random jumps \cite{Das(2017), LM18, SM10, TYM11, YM12}. The importance of this problem is linked to the fact that the dynamics of many real processes cannot be described by continuous models such as ordinary differential equations or Ito’s stochastic differential equations \cite{LM2}. More complex systems include the presence of jumps, and these jumps can occur at random $\tau_k, k\ge 1$ or deterministic time moments $t_m, m\ge 1$. In the first case, the jump-like change can be described by point processes \cite{LM3, LM6}, or in the more specific case by generalized Poisson processes, the dynamics of which is characterized only by the intensity of the jumps. The jumps of the system at deterministic moments of time $t_m$ can be described by the relation 
\begin{equation}\label{Delta}
\Delta x(t_m) = x(t_m) - x(t_m-) = g(...),
\end{equation}
where $x(t), t\ge 0$, is a random process describing the dynamics of the system. According to the works of Katz I.Ya. \cite{LM18}, Yasinsky V.K., Yurchenko I.V. and Lukashiv T.O. \cite{Lukashiv(2009)I, Lukashiv(2009)II} the description of jumps at deterministic time moments $t_m$ are quite accurately described using the relation (1). It allows a relatively simple transfer of the basic properties of stochastic systems of differential equations without jumps ($g\equiv 0$) to systems with jumps. Such properties, as will be noted below, include the Markovian property $x(t), t\ge 0$, concerning natural filtering, and the martingale properties $\|x(t)\|^2, t\ge 0$ \cite{LM1, LM5}. 
It should be noted that the description of real dynamical systems is not limited to the Wiener process and point processes (Poisson process). A more general approach is based on the use of semimartingales \cite{P04}. The disadvantage of this approach is that it is impossible to link it with the methods used for systems described by ordinary differential equations or stochastic Ito’s differential equations. The second approach to describe jump processes $x(t)$ is based on the use of semi-Markov processes, considered in the works of Korolyuk V.S. \cite{KS18} and Malyk I.V. \cite{TYM10, TYM11, TYM11_2}. In \cite{YM12}, the authors proposed the asymptotic behavior of dynamical systems described by semi-Markov processes. Several works \cite{TYM10, TYM11, TYM11_2} are devoted to the convergence of semi-Markov evolutions in the averaging scheme and diffusion approximation. The results derived in these works together with the results of the works on large deviations (e.g. \cite{KL05}) can also be used to investigate the considered problems.

Since we consider generalized classical differential equations, the approaches will also be used classical. The basic research method is based on the Lyapunov methods described in the paper by Katz I.Ya. \cite{LM18}, Lukashiv T.O., Malyk I.V. \cite{Lukashiv(2022)}. It should be noted that the application of this method makes it possible to find the optimal control for linear systems with a quadratic quality functional, which also corresponds to classical dynamical systems \cite{LM4}.

The structure of the paper is as follows. In section \ref{s_2} we consider the mathematical model of a dynamical system with jumps. It is described by a system of stochastic differential equations with Poisson's integral and external jumps. Sufficient conditions for the existence and uniqueness of the solution of this system are given there. In section \ref{s_3} we investigate the stability in probability of the solution $x(t), \ge 0$. In this section, we consider the notion of the Lyapunov function and prove the sufficient conditions for stability in probability (Theorem 1). The algorithm for computing the quality functional $J_u(y,h,x_0)$ from the known control $u(t)$ is given in the section \ref{s_4}. Moreover, we further present sufficient conditions for the existence of optimal control (Theorem 2), which are based on the existence of a Lyapunov function for the given system. Section \ref{s_5} considers constructing an optimal control for linear non-autonomous systems via the coefficients of the system. The optimal control is found by solving auxiliary Ricatti equations (Theorem 3). For the analysis of linear autonomous systems, we consider the construction of the quadratic functional. Finally, we formulate sufficient conditions of existence of an optimal control (Theorem 4), and present the explicit form such control in the case of a quadratic quality functional.


\section{Task definition}\label{s_2}

On the probability basis $(\Omega,\mathfrak{F},F,\mathbf{P})$ \cite{LM1}, consider a stochastic dynamical system of a random structure given by Ito’s stochastic differential equation (SDE) with Poisson perturbations:

$$
dx(t)=a(t-,\xi(t-),x(t-),u(t-))dt+
$$
$$
b(t-,\xi(t-),x(t-),u(t-))dw(t)+
$$

\begin{equation}\label{eq1}
\int\limits_{\mathbb{R}^m}(c(t-,\xi(t-),x(t-),u(t-),z))\widetilde{\nu}(dz,dt),~t\in \mathbb{R}_{+}\backslash K,
\end{equation}
with Markov switches
\begin{equation}\label{eq2}
\Delta x(t)\Big|_{t=t_k}=g(t_k -,\xi(t_k -),\eta_k,x(t_k -)),~~~t_k \in K=\{t_n \Uparrow\}
\end{equation}
for $\lim\limits_{n\rightarrow +\infty}t_n=+\infty$ and initial conditions

\begin{equation}\label{eq3}
x(0)=x_0\in \mathbb{R}^m,~\xi(0)=y\in\mathbf{Y}, \eta_0=h\in\mathbf{H}.
\end{equation}

Here $\xi(t), t\geq0,$ is a homogeneous continuous Markov process with a finite number of states $\mathbf{Y}:=\{y_1,...,y_N\}$ and a generator $Q$; $\{\eta_k, k\geq0\}$ is a Markov chain with values in the space $\mathbf{H}$ and the transition probability matrix $\mathbb{P}_H$; $x: [0,+\infty)\times\Omega \rightarrow\mathbb{R}^m$; $w(t)$ is an $m$-dimensional standard Wiener process; $\widetilde{\nu}(dz,dt)=\nu(dz,dt)-\mathbb{E}\nu(dz,dt)$ is a centered Poisson measure; the processes $w, \nu, \xi$ and $\eta$ are independent \cite{LM1}, \cite{LM2}. We denote by \[{ \mathfrak{F}}_{t_{k} }=\sigma (\xi(s), w(s),\nu(s,*),\eta_{e}, s\leq t_k, t_e\leq t_k) \] a minimal $\sigma $-algebra with respect to which $\xi (t)$ is measurable for all $t\in [t_{0} , t_{k} ]$ and $\eta _{n} $ for $n\le k$.

The process $x(t), t\geq0$ is $c\grave{a}dl\grave{a}g$; the control $ u(t):=u(t,x(t)):[0,T]\times\mathbb{R}^m\rightarrow\mathbb{R}^m$ is an $m$-measure function from the class of admissible controls $U$ \cite{LM4}.

The following mappings are measurable by a set of variables $a:\mathbb{R}_{+}\times\mathbf{Y}\times\mathbb{R}^m\times\mathbb{R}^m\rightarrow\mathbb{R}^m$, $b:\mathbb{R}_{+}\times\mathbf{Y}\times\mathbb{R}^m\times\mathbb{R}^m\rightarrow\mathbb{R}^m$, $c:\mathbb{R}_{+}\times\mathbf{Y}\times\mathbb{R}^m\times\mathbb{R}^m\times\mathbb{R}^m\rightarrow\mathbb{R}^m$ and function $g:\mathbb{R}_{+}\times\mathbf{Y}\times\mathbf{H}\times\mathbb{R}^m\rightarrow\mathbb{R}^m$ satisfy the Lipschitz condition

$$
|a(t,y,x_1,u)-a(t,y,x_2,u)|+|b(t,y,x_1,u)-b(t,y,x_2,u)|+
$$
$$
+\int\limits_{\mathbb{R}^m}|c(t,y,x_1,u,z)-c(t,y,x_2,u,z)|\Pi(dz)+
$$

\begin{equation}\label{eq4}
+|g(t,y,h,x_1)-g(t,y,h,x_2)|\leq L|x_1-x_2|,
\end{equation}
where $\Pi(dz)$ is defined by $\mathbb{E}\nu(dz,dt)=\Pi(dz)dt$, $L>0$, $x_1,x_2\in \mathbb{R}^m$ for $\forall t\geq0, y\in\mathbf{Y}, h\in \mathbf{H}$, and the condition

\[|a(t,y,0,u)|+|b(t,y,0,u)|+\int\limits_{\mathbb{R}^m}|c(t,y,0,u,z)|\Pi(dz)+\]
\begin{equation}\label{eq5}
|g(t,y,h,0)|\le C<\infty,
\end{equation}

The conditions defined above with respect to the mappings $a,b,c$ and $g$ guarantee the existence of a strong solution of the problem (2)--(4) with the exact stochastic equivalence \cite{Lukashiv(2022)}.

Let us denote 
\[{\it \mathbf{P}}_{k} ((y, h, x), \Gamma \times G\times \mathbf{C}):=\]
\[:=P(\xi (t_{k+1} ), \eta _{k+1}, x(t_{k+1})\in \Gamma \times G\times \mathbf{C}|(\xi (t_{k} ), \eta _{k}, x(t_k)) = (y,h,x))\] the transition probability of a Markov chain $(\xi (t_{k} ), \eta _{k} , x(t_k))$, determining the solution to the problem (2)--(4) $x(t)$, at the $k$-th step. 

\section{Stability in probability}\label{s_3}

\textbf{Definition 1.} The discrete Lyapunov operator $(lv_{k} )(y,h,x)$ on a sequence of measurable scalar functions $v_{k} (y, h, x){\it :\mathbf{Y}}\times {\it \mathbf{H}}\times {\it \mathbb{R}}^{m} \to {\it \mathbb{R}}^{1}, k\in {\it \mathbb{N}}\cup \{ 0\} $ for SDE (2) with Markov switches (3) is defined by the equation 
\[(lv_{k} )(y, h, x):= \]
\begin{equation} \label{eq6}
:= \int\limits_{{\it \mathbf{Y}}\times {\it \mathbf{H}}\times {\it \mathbb{R}}^{m} }{\it \mathbf{P}}_{k} ((y, h, x) ,du\times dz\times dl)v_{k+1} (u, z, l)-v_{k} (y, h, x), k\ge 0.
\end{equation}

When applying the second Lyapunov method to the SDE (2) with Markov switches (3), special sequences of the above mentioned functions $v_{k} (y, h, x), k\in {\it \mathbb{N}}$ are required.

\textbf{Definition 2.} The Lyapunov function for the system of the random structure (2)-(4) is a sequence of non-negative functions $\left\{v_{k} (y, h, x),k\ge 0\right\},$ for which

\begin{enumerate}
\item for all $k\ge 0, y\in {\it \mathbf{Y}}, h\in {\it \mathbf{H}}, x\in {\it \mathbb{R}}^{m} $ the discrete Lyapunov operator is defined $(lv_{k} )(y, h, x)$ (7);

\item for $r\to \infty $
$$
 \bar{v}(r)\equiv \mathop{\inf v_{k} }\limits_{\begin{array}{l} {k\in {\it \mathbb{N}}, y\in {\it \mathbf{Y}},} \\ {h\in {\it \mathbf{H}}, \left|x\right|\ge r} \end{array}} (y, h, x)\to +\infty ;
$$

\item for $r\to 0$
$$
 \underline{v}(r)\equiv \mathop{\sup v_{k} }\limits_{\begin{array}{l} {k\in {\it \mathbb{N}}, y\in {\it \mathbf{Y}},} \\ {h\in {\it \mathbf{H}}, \left|x\right|\le r} \end{array}} (y, h, x)\to 0,
$$
\end{enumerate}
and moreover $\bar{v}(r)$ and $\underline{v}(r)$ are continuous and strictly monotonous.

\textbf{Definition 3.} Let us call a system of random structure (2)-(4) stable in probability on the whole, if for $\forall \varepsilon _{1} >0, \varepsilon _{2} >0$ on can specify such $\delta>0$ that from the inequality $\left|x\right|<\delta $ follows the inequality
\begin{equation} \label{GrindEQ__1_23_}
{\it \mathbf{P}}\left\{\mathop{\sup }\limits_{t\ge 0} \left|x(t)\right|>\varepsilon _{1} \right\}<\varepsilon _{2}
\end{equation}
for all $x_0\in {\it \mathbb{R}}^{m} $, $y\in {\it \mathbf{Y}}, h\in {\it \mathbf{H}}$.

To solve the problem (2)-(4) on the intervals $[t_{k} , t_{k+1} )$ the following estimate takes place.

{\textbf{Lemma 1.}} Let the coefficients of the equation (2) $a, b, c$ and function $g$ satisfy the Lipschitz condition (5) and the uniform boundedness condition (6).

Then for all $k\ge 0$ the inequality for the strong solution of the Cauchy problem  (2)-(4) holds

\[{\it \mathbf{E}}\left\{\mathop{\sup }\limits_{t_{k} \le t< t_{k+1} } \left|x(t)\right|^{2} \right\}\le 7\left[\mathbb{E}\left|x(t_{k} )\right|^{2} +3C^{2} (t_{k+1} -t_{k} )\right]\times \]
\begin{equation} \label{eq26}
\times \exp \left\{7L^{2} ((t_{k+1} -t_{k} ) +8)\right\}, t\in(t_k, t_{k+1}).
\end{equation}

{\textit{Proof of Lemma 1}}

Using the integral form of the strong solution of the equation (2) \cite{LM5}, for all $t\in[t_k,t_{k+1}), t_k\geq0,$ the following inequality is true

$$
{\left|x(t)\right|\le \left|x(t_{k} )\right|+\int _{t_{k} }^{t}\left|a(\tau ,\xi (\tau ),x(\tau ),u(\tau ))-a(\tau ,\xi (\tau ),0,u(\tau ))\right|d\tau +}
$$
$$
+\int\limits _{t_{k} }^{t}\left|a(\tau ,\xi (\tau ),0,u(\tau ))\right|d\tau +
$$
$$
+\int\limits _{t_{k} }^{t}\left|b(\tau ,\xi (\tau ),x(\tau ),u(\tau ))-b(\tau ,\xi (\tau ),0,u(\tau ))\right|dw(\tau ) +
$$
$$
+\int _{t_{k} }^{t}\left|b(\tau ,\xi (\tau ),0,u(\tau ))\right|dw(\tau ) +
$$
$$
+\int\limits _{t_{k} }^{t}\int\limits_{\mathbb{R}^m}\left|c(\tau ,\xi (\tau ),x(\tau ),u(\tau ),z)-c(\tau ,\xi (\tau ),0,u(\tau ),z)\right|\widetilde{\nu}(dz,d\tau)+
$$
$$
+\int\limits _{t_{k} }^{t}\int\limits_{\mathbb{R}^m}\left|c(\tau ,\xi (\tau ),0,u(\tau ),z)\right|\widetilde{\nu}(dz,d\tau)
$$

Given (5), (6) and the inequality  $\left(\sum_{i=1}^n x_i\right)^2\le n\sum_{i=1}^n x_i^2$ we get:
\[\sup\limits_{t_k\leq t<t_{k+1}}\left|x(t)\right|^2\leq7\sup\limits_{t_k\leq  t<t_{k+1}}\left[  \left|x(t_k)\right|^2+\right.\]

 $$
 +\left|\int _{t_{k} }^{t}\left|a(\tau ,\xi (\tau ),x(\tau ),u(\tau ))-a(\tau ,\xi (\tau ),0,u(\tau ))\right|d\tau\right|^2+
$$
$$
+\left|\int\limits _{t_{k} }^{t}\left|a(\tau ,\xi (\tau ),0,u(\tau ))\right|d\tau\right|^2 +
$$
$$
+\left|\int\limits _{t_{k} }^{t}\left|b(\tau ,\xi (\tau ),x(\tau ),u(\tau ))-b(\tau ,\xi (\tau ),0,u(\tau ))\right|dw(\tau )\right|^2 +
$$
$$
+\left|\int _{t_{k} }^{t}\left|b(\tau ,\xi (\tau ),0,u(\tau ))\right|dw(\tau )\right|^2 +
$$
$$
+\left|\int\limits _{t_{k} }^{t}\int\limits_{\mathbb{R}^m}\left|c(\tau ,\xi (\tau ),x(\tau ),u(\tau ),z)-c(\tau ,\xi (\tau ),0,u(\tau ),z)\right|\widetilde{\nu}(dz,d\tau)\right|^2+
$$
$$
+\left.\left|\int\limits _{t_{k} }^{t}\int\limits_{\mathbb{R}^m}\left|c(\tau ,\xi (\tau ),0,u(\tau ),z)\right|\widetilde{\nu}(dz,d\tau)\right|^2\right]\leq
$$
$$
\leq7\left[\sup\limits_{t_k\leq t<t_{k+1}}\left|x(t)\right|^2+\sup\limits_{t_k\leq t<t_{k+1}}L^2\left|\int\limits _{t_{k} }^{t}\left|x(\tau)\right|d\tau\right|^2+\right.
$$
$$
+C^2(t_{k+1}-t_k)+ \sup\limits_{t_k\leq t<t_{k+1}}L^2\left|\int\limits _{t_{k} }^{t}\left|x(\tau)\right|dw(\tau)\right|^2+C^2(t_{k+1}-t_k)+
$$
$$
+\sup\limits_{t_k\leq t<t_{k+1}} \left|\int\limits _{t_{k} }^{t}\int\limits_{\mathbb{R}^m}\left|c(\tau ,\xi (\tau ),x(\tau ),u(\tau ),z)-c(\tau ,\xi (\tau ),0,u(\tau ),z)\right|\widetilde{\nu}(dz,d\tau)\right|^2+
$$
$$
+\left.\sup\limits_{t_k\leq t<t_{k+1}} \left|\int\limits _{t_{k} }^{t}\int\limits_{\mathbb{R}^m}\left|c(\tau ,\xi (\tau ),0,u(\tau ),z)\right|\widetilde{\nu}(dz,d\tau)\right|^2\right].
$$

Consider the designation
\[y(t) = \mathbb{E}\left\{\sup\limits_{t_k\leq s<t}\left|x(s)\right|^2/\mathfrak{F}_{t_k}\right\}.\]
Then, according to the last inequality, $y(t)$ satisfies the ratio  
\[y(t)\le 7\left[ \mathbb{E} \left\{\left|x(t)\right|^2/ \mathfrak{F}_{t_k} \right\} +3C^2(t_{k+1}-t_k)+L^2((t_{k+1}-t_k)+8)\cdot\int\limits _{t_{k} }^{t}y(\tau)d\tau. \right]\]

Using the Gronwall inequality, we obtain an estimate of
$$
\mathbb{E}\left\{\sup\limits_{t_k\leq t<t_{k+1}}\left|x(t)\right|^2/\mathfrak{F}_{t_k}\right\}\leq
$$
$$
\leq 7\left[\mathbb{E}\left|x(t_{k} )\right|^{2} +3C^{2} (t_{k+1} -t_{k} )\right]e^{7L^{2} \left((t_{k+1} -t_{k} ) +8\right)} ,
$$
as required to prove.

\textit{End of proof of Lemma 1.}

\textbf{Remark 1.} We will consider the stability of the trivial solution $x\equiv 0$ of the system (2)-(4) that is, the fulfillment of (6) when $C =0$ \cite{LM10}, \cite{LM8}, \cite{LM9}.

{\textbf{Theorem 1.}} Let:

1) interval lengths $[t_{k} , t_{k+1} )$ do not exceed $\Delta >0$, i.e. $\left|t_{k+1} -t_{k} \right|\le \Delta , k\ge 0;$

2) the Lipschitz condition is satisfied (5);

3) there exist Lyapunov functions $v_{k} (y, h, x), k\ge 0$ such that the following inequality holds true
\begin{equation} \label{GrindEQ__1_33_}
(lv_{k} )(y, h, x)\le 0, k\ge 0.
\end{equation}

Then the system of random structure (2)-(4) is stable in probability on the whole.

\textbf{Remark 2.} It should be noted that if condition 1 is not satisfied, the number of jumps (3) is finite and the system (2)-(4) turns into a system without jumps after $\max t_k$. In this case, we can use the results presented in \cite{LM18}.

\textit{Proof of Theorem 1.}
 The conditional expectation of the Lyapunov function is
\[{\it \mathbf{E}}\left\{{v_{k+1} (\xi (t_{k+1} ), \eta _{k+1} , x(t_{k+1} )) \mathord{\left/{\vphantom{v_{k+1} (\xi (t_{k+1} ), \eta _{k+1} , x(t_{k+1} )) {\mathfrak{F}}_{t_{k} } }}\right.\kern-\nulldelimiterspace} {\mathfrak{F}}_{t_{k} } } \right\}=\]
\begin{equation} \label{GrindEQ__1_34_}
=\int _{{\it \mathbf{Y}}\times {\it \mathbf{H}}\times {\it \mathbb{R}}^{m} }{\it \mathbf{P}}_{k}  ((\xi (t_{k} ), \eta _{k}, x(t_{k} ))(du\times dz\times dl)v_{k+1} (u, z, l)).
\end{equation}

Then by the definition of the discrete Lyapunov operator $(lv_{k} )(y, h, x)$ (see (7)) and from equation (11), taking into account (10), we can get the following inequality
\[{\it \mathbf{E}}\left\{{v_{k+1} (\xi (t_{k+1} ), \eta _{k+1} , x(t_{k+1} )) \mathord{\left/{\vphantom{v_{k+1} (\xi (t_{k+1} ), \eta _{k+1} , x(t_{k+1} )) \mathfrak{F}_{t_{k} } }}\right.\kern-\nulldelimiterspace} \mathfrak{F}_{t_{k} } } \right\}=v_{k} (\xi (t_{k} ), \eta _{k} , x(t_{k} ))+\]
\begin{equation} \label{GrindEQ__1_35_}
 +(lv_{k} )(\xi (t_{k} ), \eta _{k} , x(t_{k} ))\le \bar{v}(\left|x(t_{k} )\right|).
\end{equation}

From Lemma 1 and the properties of the function $\bar{v}$, it follows that the conditional expectation of the left part of the inequality (12) exists.

Using (11), (12), let us write the discrete Lyapunov operator $(lv_{k} )(y, h, x)$,defined on solutions (2)-(4):
\[(lv_{k}) (\xi (t_{k} ), \eta _{k} , x(t_{k} ))={\it \mathbf{E}}\left\{{v_{k+1} (\xi (t_{k+1} ), \eta _{k+1} , x(t_{k+1} )) \mathord{\left/{\vphantom{v_{k+1} (\xi (t_{k+1} ), \eta _{k+1} , x(t_{k+1} )) \mathfrak{F}_{t_{k} } }}\right.\kern-\nulldelimiterspace} \mathfrak{F}_{t_{k} } } \right\}-\]
\begin{equation} \label{GrindEQ__1_36_}
-v_{k} (\xi (t_{k} ), \eta _{k} , x(t_{k} ))\le 0.
\end{equation}

Then when $k\ge 0$ the following inequality is satisfied
\[{\it \mathbf{E}}\left\{{v_{k+1} (\xi (t_{k+1} ), \eta _{k+1} , x(t_{k+1} )) \mathord{\left/{\vphantom{v_{k+1} (\xi (t_{k+1} ), \eta _{k+1} , x(t_{k+1} )) \mathfrak{F}_{t_{k} } }}\right.\kern-\nulldelimiterspace} \mathfrak{F}_{t_{k} } } \right\}\le v_{k} (\xi (t_{k} ), \eta _{k} , x(t_{k} )).\]

This means that the sequence of random variables
$$v_{k} (\xi (t_{k} ), \eta _{k} , x(t_{k} ))$$
is a supermartingale with respect to $\mathfrak{F}_{t_{k} } $ \cite{LM6}.

Thus, the following inequality holds:
\[{\it \mathbf{E}}\left\{v_{N+1} (\xi (t_{N+1} ), \eta _{N+1} , x(t_{N+1} ))\right\} -{\it \mathbf{E}}\left\{v_{n} (\xi (t_{n} ), \eta _{n} , x(t_{n} ))\right\} =\]
\[=\sum _{k=n}^{N}{\it \mathbf{E}}\left\{(lv_{k}) (\xi (t_{k} ), \eta _{k} , x(t_{k} ))\right\}\le 0.\]

Since the random variable $\mathop{\sup }\limits_{t_{k} \le t<t_{k+1} } \left|x(t)\right|^{2} $ is independent of events of $\sigma$- algebra  $\mathfrak{F}_{t_{k} } $ \cite{LM3}, then
$$
{\it \mathbf{E}}\left\{{\mathop{\sup }\limits_{t_{k} \le t<t_{k+1} } \left|x(t)\right|^{2} \mathord{\left/{\vphantom{\mathop{\sup }\limits_{t_{k} \le t<t_{k+1} } \left|x(t)\right|^{2} \mathfrak{F}_{t_{k} } }}\right.\kern-\nulldelimiterspace} \mathfrak{F}_{t_{k} } } \right\}={ E}\left\{\mathop{\sup }\limits_{t_{k} \le t<t_{k+1} } \left|x(t)\right|^{2} \right\},
$$
i.e., the inequality (9) also holds for the usual expectation
\[{\it \mathbf{E}}\left\{\mathop{\sup }\limits_{t_{k} \le t<t_{k+1} } \left|x(t)\right|^{2} \right\}\leq 7\left[\mathbb{E}\left|x\right|^{2} \right]e^{7L^{2} \left(\Delta +8\right)} \]
at $C=0$, assuming that the stability of the trivial solution is investigated.

Then
\[ {\it \mathbf{P}}\left\{\mathop{\sup }\limits_{t\ge 0} \left|x(t)\right|>\varepsilon _{1} \right\}=\]
 \[ ={\it \mathbf{P}}\left\{\mathop{\sup }\limits_{n\in \mathbb{N}} \mathop{\sup }\limits_{t_{n-1} \le t < t_{n} } \left|x(t)\right|>\varepsilon _{1} \right\}\le \]
 \[ \le {\it \mathbf{P}}\left\{\mathop{\sup }\limits_{n\in \mathbb{N}} 7e^{7L^{2} \left(\Delta +8\right)} \left|x( t_{n-1})\right|>\varepsilon _{1} \right\}\le \]
  \[ \le {\it \mathbf{P}}\left\{\mathop{\sup }\limits_{n\in \mathbb{N}}  \left|x( t_{n-1})\right|>\frac{\varepsilon _{1}}{7} e^{-7L^{2} \left(\Delta +8\right)} \right\}\le \]
\begin{equation} \label{GrindEQ__1_39_}
\le {\it \mathbf{P}}\left\{\mathop{\sup }\limits_{n\in \mathbb{N}} v_{n-1} (\xi (t_{n-1} ),\eta _{n-1} , x(t_{n-1} )) \ge \bar{v}(\frac{\varepsilon _{1}}{7} e^{-7L^{2} \left(\Delta +8\right)} )\right\}
\end{equation}

If $\sup \left|x(t_{k} )\right|\ge r$, then based on the definition of the Lyapunov function the inequality is fulfilled
\begin{equation} \label{GrindEQ__1_40_}
 \mathop{\sup }\limits_{k\ge 0} v_{k} (\xi (t_{k} ),\eta _{k} , x(t_{k} ))\ge \inf\limits_{k\ge 0, y\in {\it \mathbf{Y}}, {h\in {\it \mathbf{H}}, \left|x\right|\ge r} } v_{k} (y, h,x) =\bar{v}(r) .
\end{equation}

Using the inequality for non-negative supermartingales \cite{LM1}, \cite{LM6}, we obtain an estimate of the right-hand side (14):
\[{\it \mathbf{P}}\left\{\mathop{\sup }\limits_{n\in \mathbb{N}} v_{n-1} (\xi (t_{n-1} ),\eta _{n-1} , x(t_{n-1} ))\ge \bar{v}(\frac{\varepsilon _{1}}{7} e^{-7L^{2} \left(\Delta +8\right)})\right\}\le \]
\begin{equation} \label{GrindEQ__1_41_}
\le \frac{1}{\bar{v}(\frac{\varepsilon _{1}}{7} e^{-7L^{2} \left(\Delta +8\right)} )} v_{0} (y, h, x)\le \frac{\bar{v}(\left|x\right|)}{\bar{v}(\frac{\varepsilon _{1}}{7} e^{-7L^{2} \left(\Delta +8\right)})} .
\end{equation}

Given the inequality (14), the inequality (16) makes it possible to guarantee that the inequality (8) of stability in probability for the whole system (2)-(4).

\textit{End of proof of Theorem 1.}

\section{Stabilization}\label{s_4}

The optimal stabilization problem is that for an SDE (2) with switches (3), one should construct such a control $u(t,x(t))$ that the unperturbed motion $x(t)\equiv0$ of the system (2)-(4) is stable in probability on the whole.

It is assumed that the control $u$ will be determined by the full feedback principle. In addition, the condition of continuity of $u(t)$ on $t$ in the range

\begin{equation}\label{eq17}
t\geq 0,~ x\in \mathbb{R}^m,~y\in\mathbf{Y}, ~h\in\mathbf{H}.
\end{equation}
for every fixed $\xi(t)=y\in\mathbf{Y}$ and $\eta_{k}=h\in\mathbf{H}$.

It is also assumed that the structure in which the system is at time $t\ge 0$, which is independent of the Markov chain $\eta_k$ ($k\ge 0$ corresponds to time $t_k\in K$), is known.

Obviously, there is an infinite set of controls. The only control should be chosen from the requirement of the best quality of the process, which is expressed in the form of the minimization condition of the functional

\[{\it I}_{u} (y,h,x_0):= \]
\begin{equation} \label{eq18}
:= \sum _{k=0}^{\infty }\int _{t_k}^{\infty }{\it E}  \left\{ W(t, x(t) ,u(t)) / \xi(0)=y, \eta_0=h,x(0)=x_0 \right\} dt,
\end{equation}
where $W(t,x,u)\ge 0$ -- is a non-negative function defined in the region $t\ge 0, x\in {\it \mathbb{R}}^{m}$, $u\in {\it \mathbb{R}}^{r} $.

The algorithm for calculating the (18) functional for a given control $u(t,x)$ is as follows.

A) Find the trajectory $x(t)$ with an SDE (2) at $u\equiv u(t,y,h,x)$, for example by the Euler-Maruyama method \cite{KP92}.

B) Substitute $x(t)$, $\xi (t)$, $u(t)=u(t,x(t))$ into the functional (18).

C) Calculate the value of the function (18) by statistical modeling (Monte Carlo).

D) The problem of choosing the functional $W(t,x,u)$, which determines the estimate ${\it I}_{u} $ and the quality of the process $x(t)$ as a strong solution of the SDE (2), is related to the specific features of the problem and the following three conditions can be identified:

\begin{enumerate}
\item   the minimization conditions of the functional (18) must ensure that the strong solution $x(t)$ of the SDE (2) fades fast enough on average with high probability;
\item	the value of the integral should satisfactorily estimate the computation time spent on generating the control $u(t)$;
\item	the value of the quality functional should satisfactorily estimate the computation time spent on forming the control $u(t)$;
\item	The functional $W(t,x,u)$ must be such that the solution of the stabilization problem can be constructed.
\end{enumerate}

\textbf{Remark 3.}
For linear SDE (2), in many cases the quadratic form with respect to the variables $x, u$ is satisfactory
\begin{equation}\label{eq_control}
W(t,x,u)=x^{T} M(t)x+u^{T} D(t)u, 
\end{equation}
where $M(t)$ -- is a symmetric non-negative matrix of size $m\times m$; $D(t)$ -- is a positively determined matrix of the size $r\times r$ for all $t\ge 0$.

\textbf{Remark 4.}
Note that according to the feedback principle, $M(t)$ and $D(t)$ depend indirectly on the values of $\xi(t)$ and $\eta_k$. Therefore, in the examples below, we will calculate the values of $M(t)$ and $D(t)$ for fixed $\xi(t)$ and $\eta_k$.

The value ${\it I}_{u} $ in the case of the quadratic form of the variables $x$ and $u$ evaluates the quality of the transition process quite well on average. The presence of the term $u^{T} Du$ and the minimum condition simultaneously limit the amount of the control action $u\in {\it \mathbb{R}}^{r} $.

\textbf{Remark 5.}
If the jump condition of the phase trajectory is linear, then the solution of the stabilization problem belongs to the class of linear on the phase vector $x\in {\it \mathbb{R}}^{m} $ controls $u(t,x)$. Such problems are called linear-quadratic stabilization problems.

\textbf{Definition 4.}
The control $u^0(t)$, which satisfies the condition 
$$
{\it I}_{u^0}(y,h,x_0)= \min {\it I}_u(y,h,x_0),
$$
where the minimum should be searched for all controls continuous variables $t$ and $x$ at each $\xi(0)=y\in {\it \mathbf{Y},}  \eta_0=h\in {\it \mathbf{H}}$, let us call it optimal in the sense of optimal stabilization of the strong solution $x\in {\it \mathbb{R}}^m$ of the system (2)-(4).

\textbf{Theorem 2.}
Let the system (2)-(4) have a scalar function $v^0(t_k,y,h,x)$ and $r$-vector function $u^0(t,y,h,x)\in {\it \mathbb{R}}^r $ in the region (17) and fulfill the conditions:

\begin{enumerate}
\item	the sequence of the functions $v_k^{0} (y,h,x)\equiv v^{0} (t_k,y,h,x)$ is the Lyapunov functional; 
\item	the sequence of $r$-measured functions-control
\begin{equation} \label{eq20}
u_{k}^{0} (y,h,x)\equiv u^0 (t_k,y,h,x)\in \it \mathbb{R}^{r};
\end{equation}
is measurable in all arguments, where $0\le t_{k} <t_{k+1} ,    k\ge 0$;
\item   sequence of functions from the criterion (18) by $x\in {\it \mathbb{R}}^{m}$ is positive definite, i.e. for $\forall t\in [t_{k} ,t_{k+1} ),    k\ge 0$,
\begin{equation} \label{eq21}
W(t,x,u_{k}^{0} (y,h,x))>0;
\end{equation}
\item   the sequence of infinitesimal operators $(lv_{k}^{0} )\left. \right|_{u_{k}^{0} } $, calculated for $u_{^{k} }^{0} \equiv u^{0} (y,h,x)$, satisfies the condition for $\forall t\in [t_{k} ,t_{k+1} )$
\begin{equation} \label{eq22}
(lv_{k}^{0} )\left. \right|_{u_{k}^{0} } =-W(t,x,u_{k}^{0} );
\end{equation}
\item   the value of $(lv_{k}^{0} )+W(t,x,u)$ reaches a minimum at $u=u^{0}, k\ge 0$, i.e.
\[(lv_{k}^{0} )\left. \right|_{u_{k}^{0} } +W(t,x,u_{k}^{0} )=\]
\begin{equation} \label{eq23}
=\mathop{\min }\limits_{u\in {\it \mathbb{R}}^{r} } \{ (lv_{k}^{0} )\left. \right|_{u} +W(t,x,u)\} =0.
\end{equation}
\item   the series 
\begin{equation} \label{eq24}
\sum _{k=0}^{\infty }\int _{t_{k} }^{\infty }{\it \mathbf{E}}  \left\{W(t,x(t),u(t))/ x(t_{k-1} ) \right\}dt<\infty
\end{equation}
converges.
\end{enumerate}

Then the control $u_{^{k} }^{0} \equiv u^{0} (t_{k},y,h,x),    k\ge 0$, stabilizes the solution of the problem (2)-(4). In this case, is hold the equality 
\[v^{0} (y,h,x_0)\equiv \]
\[\equiv \sum _{k=0}^{\infty }\int _{t_{k} }^{\infty }{\it \mathbf{E}}  \left\{W(t,x(t),u(t))/ x(t_{k-1} ) \right\}dt=\]
\begin{equation} \label{eq25}
=\mathop{\min }\limits_{u\in {\it \mathbb{R}}^{r} } \sum _{k=0}^{\infty }\int _{t_{k} }^{\infty }{\it \mathbf{E}}  \left\{W(t,x(t) ,u(t))/x_(t_{k})\right\}dt\equiv {\it I}_{u^{0} } (y,h,x_0). 
\end{equation}

\textit{Proof of Theorem 2.}

\textbf{I.} Stability in probability in the whole of a dynamical system of a random structure (2)-(4) for  $u\equiv u^{0} (t_k,x),  k\ge 0$ immediately follows from the Theorem 1, since the functionals $v^{0} (y,h,x)$ for anyone $t\in [t_{k} ,t_{k+1}),    k\ge 0$ satisfy the conditions of this theorem.

\textbf{II.} The equality (25) is obviously also a consequence of the Theorem 1.

\textbf{III.} Proof by contradiction that the stabilization of a strong solution of a dynamical system of random structure (2)-(4) is controlled by $u^{0} (t_k,x),$ $t_{k} \le t<t_{k+1} ,k\ge 0$.

Let exist the control $u^{*} (t_{k},x)\ne u^{0} (t_{k},x)$, which, when substituted into the SDE (2), realizes a solution $x^{*} (t)$ with initial conditions (3), (4) such that hold the equality
\begin{equation} \label{eq26}
{\it I}_{u^{*} } (y,h,x_0)\le {\it I}_{u^{0}} (y,h,x_0).
\end{equation}

The fulfillment of conditions 1)-6) of the Theorem 2 will lead to an inequality (see (27)) opposite to (26). 

From the condition 5) (see (23)) follows the inequality
\begin{equation} \label{eq27}
(lv_{k}^{0} )\left. \right|_{u^{*} } \ge -W(t,x,u^{*} (t,y,h,x)).
\end{equation}

Averaging (27) over random variables $\{ x^{*} (t),\xi (t),\eta _{k} \} $ over intervals $[t_{k} ,t_{k+1}),  k\ge 0$ and integrating over $t$ from $0 $ to $T$, we obtain $n$ inequalities:
\[{\it \mathbf{E}}{  \{ }v^{0} {  (}{t_{1} ,\xi (t_{1} ),\eta _{k_{1} },  x^{*} (t_{1} )) \mathord{\left/{\vphantom{t_{1} ,\xi (t_{1} ), \eta _{k_{1} }, x^{*} (t_{1} ), ) y ,h,x_0 }}\right.\kern-\nulldelimiterspace} y_{1} ,\eta _{k_{1} }, x^{*} (t_{1}) } {  \} -} v^{0}(y ,h,x_0)\ge \]
\begin{equation} \label{eq28}
\ge -\int _{t_{0} }^{t_{1} }{\it \mathbf{E}} \left\{W(t,x^{*} (t),u^{*} (t))/x_0 \right\} dt,       
\end{equation}
\[{\it \mathbf{E}}{  \{ }v^{0} {  (}{t_{2} ,\xi (t_{2} ),\eta _{k_{2}}, x^{*} (t_{2} )  ) \mathord{\left/{\vphantom{t_{2} ,\xi (t_{2} ),  \eta _{k_{2} }, x^{*} (t_{2} ) ) y_{1} ,\eta _{k_{1} },x^{*} (t_{1}) }}\right.\kern-\nulldelimiterspace} y_{1} ,\eta _{k_{1} },x^{*} (t_{1}) } {  \} -}\]
$$
{  -\{ }v^{0} {  (}{t_{1} ,\xi (t_{1} ),  \eta _{k_{1} }, x^{*} (t_{1}) ) \mathord{\left/{\vphantom{t_{1} ,\xi (t_{1} ), \eta _{k_{1} }, x^{*} (t_{1} ) ) y ,h, x_0 }}\right.\kern-\nulldelimiterspace} y,h,x_0 } {  \} }\ge 
$$
\begin{equation} \label{eq29}
\ge -\int _{t_{1} }^{t_{2} }{\it \mathbf{E}} \left\{W(t,x^{*} (t),u^{*} (t))/x^{*} (t_{1}) \right\} dt, 
\end{equation}

\begin{center}
...    
\end{center}

\[{\it \mathbf{E}}{  \{ }v^{0} {  (}{t_{n} ,\xi (t_{n} ), \eta _{k_{n} }, x^{*} (t_{n}) ) \mathord{\left/{\vphantom{t_{n} ,\xi (t_{n} ),  \eta _{k_{n} },x^{*} (t_{n}) ) y_{n-1} ,x^{*} (t_{n-1}),\eta _{k_{n-1} } }}\right.\kern-\nulldelimiterspace} y_{n-1},\eta _{k_{n-1} }, x^{*} (t_{n-1}) } {  \} -}\]
\[{  -\{ }v^{0} {  (}{t_{n-1} ,\xi (t_{n-1} ), \eta _{k_{n-1} },x^{*} (t_{n-1}) ) \mathord{\left/{\vphantom{t_{n-1} ,\xi (t_{n-1} ), \eta _{k_{n-1} },x^{*} (t_{n-1}) ) y_{n-2} ,\eta _{k_{n-2} },x^{*} (t_{n-2}) }}\right.\kern-\nulldelimiterspace} y_{n-2} ,\eta _{k_{n-2} }, x^{*} (t_{n-2}) } {  \} }\ge \]
\begin{equation}\label{eq30}
\ge -\int _{t_{n-1} }^{t_{n} }{\it \mathbf{E}} \left\{W(t,x^{*} (t),u^{*} (t))/x^{*} (t_{n-1})  \right\} 
\end{equation}
 
Taking into account the martingale property of the Lyapunov functions $v^{0} (t,\xi (t),h,x^{*} (t))$ (see condition 1) of the theorem) due to the system (2)-(4), i.e. by the definition of a martingale, we have $n$ equalities with the probability of one:
\[{\it \mathbf{E}} \{ v^0({t_{k} ,\xi (t_{k} ), \eta _{k}, x^{*} (t_{k})) / y_{k-1} ,\eta _{k-1},x^{*} (t_{k-1} ) \} }=\]
\begin{equation} \label{eq31}
=v^{0} (t_{k-1} ,\xi (t_{k-1} ),  \eta _{k-1},x^{*} (t_{k-1}) ),                        k=\overline{1,n}.
\end{equation}

Substituting (31) into the inequalities (28)-(30), we obtain the inequality
$$
{\it \mathbf{E}}{  \{ }v^{0}(t_{n} ,\xi (t_{n} ),  \eta _{k_{n} }, x^{*} (t_{n}) ) /{t_{n-1} ,\xi (t_{n-1} ), \eta _{k_{n-1} } ,x^{*} (t_{n-1})  \} -}v^{0} {  (}y,h,x_0)\ge
$$
$$
\ge -\sum _{k=0}^{n}\int _{t_{k} }^{t_{k+1} }{  \mathbf{E}}  \left\{W(t,  x^{*} (t),u^{*} (t)) /x^{*} (t_{k-1})  \right\}dt\ge
$$

\begin{equation} \label{eq32}
\ge -\sum _{k=0}^{\infty }\int _{t_{k} }^{\infty }{\it \mathbf{E}}  \left\{W(t,x^{*} (t),u^{*} (t)) / x^{*} (t_{k-1}) \right\}dt.
\end{equation}

According to the assumption (26) it follows that for $t_{n} \to \infty $ the integrals on the right-hand side (32) converge and, taking into account the convergence of the series (24) (condition 6)), we have the inequality:
$$
v^{0} {  (} y,h,x_0 )={\it I}_{u^{0} } (y,h,x_0)\le
$$
$$
\le \sum _{k=0}^{\infty }\int _{t_{k} }^{\infty }{\it \mathbf{E}}  \left\{W(t, x^{*} (t),u^{*} (t)) /x^{*} (t_{k-1}) \right\}dt=
$$
\begin{equation} \label{eq33}
={\it I}_{u^{*} } (y,h,x_0).
\end{equation}

Indeed, from the convergence of the series (32), under condition 6) it follows that the integrands in (33) tend to zero as $t\to \infty $. In this way, $\mathop{\lim }\limits_{n\to \infty } {\it \mathbf{E}}{  \{ }v^{0} {  (}t_{n} ,y_{n},\eta _{k_{n} },x^{*} (t_{n} ) {  \} }={  0}$.

Note that it makes sense to consider natural cases when from the condition 
$$
{\it \mathbf{E}}{  \{ }W{  \} }\mathop{\to }\limits_{{t}\to \infty } {0}
$$
it follows that ${\it \mathbf{E}}{\{ }v^{{0}} {  \} }\mathop{\to }\limits_{{t}\to \infty } {0}$.

Thus, the inequality (33) contradicts the inequality (26). This contradiction proves the statement about the optimality of the control  $u^{0} (t_{k} ,x),    k\ge 0$. 

\textit{End of proof of Theorem 2.}

In case when the Markov process with a finite number of states $\xi (t_{k} )$ admits a conditional expansion of the conditional transition probability

\[  {\it \mathbf{P}}\{ \omega :\xi (t+\Delta t)=y_{j} /\xi (t)=y_{i} ,  y_{i} \ne y_{j} \} =\]
\begin{equation} \label{eq34}
           =q_{ij} (t)\Delta t+o (\Delta t),            i,j=\overline{1,N},
\end{equation}
we obtain an equation that must be satisfied by the optimal Lyapunov functions $v_{k}^{0} (y,h,x)$ and the optimal control $u_{k}^{0} (t,x),    \forall t\in [t_{k} ,t_{k+1} )$.

Note that according to \cite{Lukashiv(2016)}, \cite{Lukashiv(2009)I}, \cite{Lukashiv(2009)II}  the weak infinitesimal operator (7) has the form

$$
(lv_k)(y,h,x) = \frac{\partial v_k (y,h,x)}{\partial t}+(\nabla v_k (y,h,x),a(t,y,x,u))+
$$
$$
+\frac{1}{2}Sp(b^T(t,y,x,u)\cdot \nabla^2 v_k (y,h,x)\cdot b(t,y,x,u))+
$$
$$
+\int\limits_{{\it \mathbb{R}}^{m}}[v_k (y,h,x+c(t,y,x,u,z))-v_k (y,h,x)-(\nabla v_k (y,h,x))^T\cdot c(t,y,x,u,z)]\Pi(dz)+
$$
\begin{equation}\label{eq35}
+\sum\limits_{j\neq i}^{N}[\int\limits_{{\it \mathbb{R}}^{m}} v_j (t,x)p_{ij}(t,z/x)dz-v_i (t,x)]q_{ij},    
\end{equation}
$(\cdot,\cdot)$ is a scalar product, $\nabla v_k=\left( \frac{\partial v_k}{\partial x_1},...,\frac{\partial v_k}{\partial x_m} \right)^T$, $\nabla^2 v_k=\left[\frac{\partial^2 v_k}{\partial x_i \partial x_j}\right]_{i,j=1}^{m},k\geq 0$, ''$T$'' stands for a transposition, $Sp$ is a trace of matrix, $p_{ij}(t,z/x)$ is a conditional probability density
$$
P{x(\tau)\in[z,z+dz]/x(\tau-0)=x}=p_{ij}(\tau,z/x)dz+o(dz)
$$
assuming $\xi(\tau-0)=y_i, ~\xi(\tau)=y_j$.

Taking into account the formula (35), the first equation for $v^{0}$ can be obtained by putting the expression for the averaged infinitesimal operator $\left. (l v_{k}^{0} )\right|_{u^{*} } $ \cite{LM18} into the left side of (23).

Then the desired equation at the points $(t_{k} ,y_{j}, \eta_k ,x)$ has the form
\[  \frac{\partial v_{k}^{0} }{\partial t} +\left(\left(\frac{\partial v_{k}^{0} }{\partial x} \right)^{T} \cdot a(t,y,x,u)\right)+\frac{1}{2} Sp\left( \left(b^T (t,  y_{i} ,  x)\cdot \frac{\partial ^{2} v_{k}^{0}}{\partial x^{2} }\cdot b (t,  y_{i} ,  x)\right)\right)+\]
$$
+\int\limits_{{\it \mathbb{R}}^{m}}[v_k^0 (\cdot,\cdot,x+c(t,y,x,u,z))-v_k^0 -(\frac{\partial v_{k}^{0}}{\partial x})^T\cdot c(t,y,x,u,z)]\Pi(dz)+
$$
 \[+    \sum _{j\ne i}^{l}\left[\int _{-\infty }^{+\infty }v_j^{0} (y_{j}, h, x_{j}) p_{ij} \left(t,  z/x\right)dz-v_i^{0} (y_{i} ,h,x)\right]q_{ij} \left(t\right)dt +\]
\begin{equation} \label{eq36}
+W(t,x,u)=0.
\end{equation}

The second equation for optimal control  $u_{k}^{0} (t,y,h,x)$ we obtain from (36) by differentiation with respect to the variable $u$, since $u=u^{0} $ delivers the minimum of the left side of (36)
\begin{equation} \label{eq37}
\left. \left[\left(\frac{\partial v^{0} }{\partial x} \right)^{T} \cdot \left(\frac{\partial a}{\partial u} \right)+\left(\frac{\partial W}{\partial u} \right)^{T} \right]\right|_{u=u_{k}^{0} } =0,
\end{equation}
where $\frac{\partial a}{\partial u} $ -- $m\times r$-matrix of  Jacobi, stacked with elements $\left\{\frac{\partial a_{n} }{\partial u_{s} } ,  n=\overline{1,m  },  s=\overline{1,r}\right\};$$\left(\frac{\partial W}{\partial u} \right)\equiv \left(\frac{\partial W}{\partial u_{1} } ,...,\frac{\partial W}{\partial u_{r} } \right),    k\ge 0$.

Thus, the problem of optimal stabilization, according to the Theorem 2, consists in solving of a complex nonlinear system of equations (23) with partial derivatives to determine the unknown Lyapunov functions $v_{ik}^{0} \equiv v_{k}^{0} (y,h,x),    i=\overline{1,l},  k\ge 0 $.

Note that this system is obtained by eliminating the control  $u_{k}^{0} =u^{0} (t,y,h,x)$ from the equation (36), (37).

It is quite difficult to solve such a system, therefore, we will further consider linear stochastic systems for which convenient solution schemes can be constructed.

\section{Stabilization of linear systems}\label{s_5}

Consider a controlled stochastic system defined by a linear Ito's SDE with Markov parameters and Poisson perturbations 
$$
dx(t)=[A(t-,\xi (t-))x(t-)+B(t-,\xi (t-))u(t-)]dt+\sigma (t-,\xi (t-))x(t-)dw(t) + 
$$
\begin{equation}\label{eq38}
\int\limits_{\mathbb{R}^m}c(t-,\xi(t-),u(t-),z)x(t-)\widetilde{\nu}(dz,dt),~t\in \mathbb{R}_{+}\backslash K,
\end{equation}
with Markov switching
\begin{equation}\label{eq39}
\Delta x(t)\Big|_{t=t_k}=g(t_k -,\xi(t_k -),\eta_k,x(t_k -)),~~~t_k \in K=\{t_n \Uparrow\}
\end{equation}
for $\lim\limits_{n\rightarrow +\infty}t_n=+\infty$ and initial conditions

\begin{equation}\label{eq40}
x(0)=x_0\in \mathbb{R}^m,~\xi(0)=y\in\mathbf{Y}, \eta_{0}=h\in\mathbf{H}.
\end{equation}

Here $A,B,\sigma,C$ are piecewise continuous integrable matrix functions of appropriate dimensions.

Let us assume that the conditions for the jump of the phase vector $x\in {\it \mathbb{R}}^{m} $  at the moment $t=t^{*} $ of the change in the structure of the system due to the transition $\xi (t^{*} -)=y_{i} $ in $\xi (t^{*} )=y_{j} \ne y_{i} $ are linear and given in the form
\begin{equation} \label{eq41}
x(t^{*} )=K_{ij} x(t^{*} -)+\sum _{s=1}^{N}\xi _{s} Q_{s}  x(t^{*} -),
\end{equation}
where $\xi _{s} := \xi _{s} (\omega )$ -- independent random variables for which  ${\it \mathbf{E}}\xi _{s} =0,  {\it \mathbf{E}}\xi _{s}^{2} =1$,  $K_{ij} $ and $Q_s$  -- given $(m\times m)$-matrices.

Note that the equality (41) can replace the general jump conditions \cite{Lukashiv(2009)I}:

- the case of non-random jumps will be at $Q_{s} =0$, i.e.
\[x(t^{*} )=K_{ij} x(t^{*} -);\]

- continuous change of the phase vector means that $Q_{s} =0$, $K_{ij} =A_{ij} =I$ (identity $(m\times m)$-matrix).

The quality of the transition process will be estimated by the quadratic functional
\begin{equation} \label{eq42}
{\it I}_{u} (y,h,x_0):=
\end{equation}
$$
:= \sum _{k=0}^{\infty }\int _{t_{k} }^{\infty }{\it \mathbf{E}}\left\{x^{T} (t) M(t)x(t)+u^{T} (t)D(t)u(t)\right.   /\left. y,h,x_0 \right\}dt,
$$
where $M(t)\ge 0,    D(t)>0$ -- symmetric matrices of dimensions $(m\times m)$ and $(r\times r)$ respectively.

According to the Theorem 2 we need to find optimal Lyapunov functions $v_{k}^{0} (y,h,x)$ and a control $u_{k}^{0} (t,x)$   for $    \forall   t\in [t_{k} ,t_{k+1} ),    t_{k} \in K,    k=0,1,2,...$.

The optimal Lyapunov functions are sought in the form
\begin{equation} \label{eq43}
v_{k}^{0} (y,h,x)=x^{T} G(t,y,h)x,
\end{equation}
where $G(t,y,h)$ -- positive-definite symmetric matrix of the size $(m\times m)$.

Everywhere below, when $\xi (t)$ describes a Markov chain with a finite number of states ${\it \mathbf{Y}}= \{ y_{1} ,  y_{2} ,  ...,  y_{l} \} $, and $\eta _{k} ,  k\ge 0$, -- a Markov chain with  values $h_{k}$ in metric space ${\it \mathbf{H}}$ and with transition probability at the $k$-th step ${\it \mathbf{P}}_{k} (h,  G)$, we introduce the following notation:
$$ 
A_{i}(t):= A(t,y_{i}), ~   B_{i}(t):= B(t,y_{i}),~ \sigma_i (t):=\sigma(t,y_i), ~C_i(t,z):=C(t,y_i,z),
$$
$$ 
G_{ik} (t):= G(t,y_{i} ,h_{k} ),~~v_{ik} := v(y_{i} ,h_{k},x ).
$$

Let us substitute the functional (43) in equations (36) and (38) to find the an optimal Lyapunov function $v_{k}^{0} (y,h,x)$ and an optimal control $u_{k}^{0} (t,x)$ for $  \forall   t\in [t_{k} ,t_{k+1} )$ and, given the form of a weak infinitesimal operator (35), we obtain:
$$
x^T(t)\frac{dG_{ik}(t)}{dt}x(t)+2[A_i(t)x(t)+B_i(t)u(t)]G_{ik}(t)x(t)+
$$
$$
+Sp(x^T(t)\sigma_i^T(t)G_{ik}(t)\sigma_i(t)x(t))+\int\limits_{\mathbb{R}^m}x^T(t)C_i^T(t,z)G_{ik}(t)C_i(t,z)x(t)\Pi(dz)+
$$
$$
+x^T(t)\sum\limits_{j\neq i}^{N}\left[K_{ij}^T G_{ik}(t)K_{ij}+\sum\limits_{s=0}^{l}Q_s^T G_{ik}(t)Q_s-G_{ik}(t)\right]q_{ij}x(t)+
$$
\begin{equation}\label{eq44}
    +x^T(t)M_{ik}(t)x(t)+u^T(t)D_{ik}(t)u(t)=0,
\end{equation}

\begin{equation}\label{eq45}
  2x^T(t)G_{ik}(t)B_i(t)+2u^T(t)D_{ik}(t)=0.
\end{equation}

Note that the partial derivative with respect to $u$ of the operator $(lv)$ is equal to zero, which confirms the conjecture about constructing an optimal control that does not depend on switching (39) for the system (40).

From (45) we find an optimal control for $\xi (t)=y_{i} $, when switching (39) $\eta _{k} =h_{k} ,  k\ge 0,$
\begin{equation} \label{eq46}
u_{ik}^{0} (t,x)=-D_{ik}^{-1}(t) B_{i}^{T}(t) G_{ik}(t) x(t).
\end{equation}

Given the matrix equality
$$
2x^{T}(t) G_{ik}(t) A_{i}(t) x=x^{T}(t) (G_{ik}(t) A_{i}(t) +A_{i}^{T}(t) G_{ik}(t) )x(t),
$$
excluding $u_{ik}^{0}$ from (44) and equating to zero the resulting matrix of a quadratic form, we can obtain a system of matrix differential equations of Riccati type for finding the matrices $G_{ik}(t) $, where $i=1,2,...,l$,  $k\ge 0$, corresponds to the interval $[t_{k} ,t_{k+1} )$:

$$
\frac{dG_{ik}(t)}{dt}+G_{ik}(t)A_i(t)-B_i(t)D_{ik}^{-1}(t) B_{i}^{T}(t) G_{ik}(t))+
$$
$$
+Sp(\sigma_i^T(t)G_{ik}(t)\sigma_i(t))+\int\limits_{\mathbb{R}^m}C_i^T(t,z)G_{ik}(t)C_i(t,z)\Pi(dz)+
$$
\begin{equation}\label{eq47}
+\sum\limits_{j\neq i}^{N}\left[K_{ij}^T G_{ik}(t)K_{ij}+\sum\limits_{s=0}^{l}Q_s^T G_{ik}(t)Q_s-G_{ik}(t)\right]q_{ij}+M_{ik}(t)=0,
\end{equation}

\begin{equation}\label{eq48}
\lim_{t\to\infty}G_{ik}(t)=0, i=\overline{1,N}, k\geq 0.
\end{equation}

Thus, we have received the following statement.

\textbf{Theorem 3.}
Let the system of matrix equations  (47), (48) has positive-definite solutions of the order $(m\times m)$
$$
G_{1k}(t) >0,  G_{2k}(t) >0,...,G_{lk}(t) >0.
$$

Then the control (46) gives a solution to the problem of optimal stabilization of the system (38)-(40) with jump condition (41) and the criterion of optimality (42).

\section{Stabilization of autonomous systems}\label{s_6}

Consider the case of an autonomous system that is given by the SDE
\begin{equation}\label{eq49}
dx=[A(\xi (t))x+B(\xi (t))u]dt+\sigma (\xi (t))xdw(t) + C(\xi(t))xdN(t),~t\in \mathbb{R}_{+}\backslash K,
\end{equation}
with Markov switching (39) and initial conditions (40). Here $ x\in   {\it \mathbb{R}}^{m} $,  $u\in {\it \mathbb{R}}^{r} $,  $A(y)$, $ B(y)$,  $\sigma (y)$, $C(y)$ are known matrix functions defined on the set ${\it \mathbf{Y}}= \{ y_{1} ,  y_{2} ,  ...,  y_{k} \} $ of possible values of the Markov chain $\xi$. $N(t), t\geq 0$ -- Poisson process with intensity $\lambda$ \cite{LM3}.

In the case of phase vector jumps (41) and the quadratic quality functional (42) the system (47), (48) for finding unknown matrices $G_{ik}, i=\overline{1,N}, k\geq 0,$ will take the form

$$
G_{ik}A_i+A_i^T G_{ik}-B_iD_{ik}^{-1} B_{i}^{T} G_{ik}+\sigma_i^TG_{ik}\sigma_i+
$$
$$
+\lambda C_i^TG_{ik}C_i+
$$
\begin{equation}\label{eq50}
+\sum\limits_{j\neq i}^{N}\left[K_{ij}^T G_{ik}K_{ij}+\sum\limits_{s=0}^{l}Q_s^T G_{ik}Q_s-G_{ik}\right]q_{ij}+M_{ik}=0, i=\overline{1,N}, k\geq 0.
\end{equation}

\textbf{Remark 6.}
Note that any differential system written in the normal form (such as the system (38), where the dependence of $x$ on $t$ is explicitly indicated) can be reduced to an autonomous system by increasing the number of unknown functions (coordinates) by one.

\subsection{Small parameter method for solving the problem of the optimal stabilization}

The possibility of algorithmic solution of the problem of optimal stabilization of a linear autonomous system of random structure  (43), (39) (40) is achieved by introducing a small parameter \cite{LM18}. There are two ways to introduce the small parameter:

\textbf{Case I.} Transition probabilities $y_{i} \to y_{j} $ of Markov chains $\xi$ are small, i.e. the transition intensities $q_{ij}$ due to the small parameter $\varepsilon >0$ can be represented as
\begin{equation} \label{eq51}
q_{ij} =\varepsilon r_{ij} .
\end{equation}

\textbf{Case II.} Small jumps of the phase vector $x(t)\in {\it R}^{m} $, i.e. matrices $K_{ij} $ and $Q_{s} $ from (41), should be presented in the form
\begin{equation} \label{eq52}
K_{ij} =I+\varepsilon K_{ij} ; ~  Q_{s} =\varepsilon   Q_{s} .
\end{equation}

In these cases, we will search the optimal Lyapunov function $v_{k}^{0} (y,x,h), k\ge 0$, in the form of a convergent power series with a base $\varepsilon >0$
\begin{equation} \label{eq53}
v_{k}^{0} (y,h,x)=x^{T} \sum _{r=0}^{\infty }\varepsilon ^{r} G^{(r)} (y,h)x .
\end{equation}

According to (46), the optimal control $u^{0}$ should be sought in the form of a convergent series
\begin{equation} \label{eq54}
u_{k}^{0} (y,h,x)=-[D^{-1} (y)B^{T} (y)\sum _{r=0}^{\infty }\varepsilon ^{r} G^{(r)} (y,h)]x .
\end{equation}

\textbf{Case I. } Let us substitute the series (53), (54), taking into account (51), into (44):
$$
G_{ik} A_i +(A_i )^{T} G_{ik} -B_i D_{ik}^{-1} B_i^{T} G_{ik} +\sigma _i^{T} G_{ik} \sigma _i +
$$
$$
+\lambda C_i^TG_{ik}C_i+
$$
$$
+\sum \limits_{j\ne i}^{l}K_{ij}^{T} G_{ik} K_{ij}  +\sum _{s=1}^{N}Q_{s}^{T} (G_{ik} Q_{s} -G_{ik} )\varepsilon r_{ij}  +C_{ik} =0;    i=\overline{1,l},    k\ge 0.
$$

Equating the coefficients at the same powers of $\varepsilon>0$, we get:
\[A_i^{T} G_{i}^{(0)} +G_{ik}^{(0)} A_i - B_i D_{ik}^{-1} B_{ik}^{T} G_{ik}^{(0)} +\]
\begin{equation} \label{eq55}
+\sigma _{ik}^{T} G_{ik}^{(0)} \sigma_{ik}+ \lambda C_i^TG_{ik}^{(0)}C_i + M_{ik} =0,    i=\overline{1,l},    k\geq0,
\end{equation}
\[\begin{array}{l} {                                          \tilde{A}_{ik}^{T} G_{ik}^{(r)} +G_{ik}^{(r)} \tilde{A}_{ik} +\sigma _{i}^{T} G_{ik}^{(r)} \sigma _{i} + \lambda C_i^TG_{ik}^{(0)}C_i =} \\ {=-\sum _{j\ne i}^{l}(K_{ik}^{T} G_{ik}^{(r-1)} K_{ij} +\sum _{s=1}^{N}Q_{s}^{T} G_{ik}^{(r-1)} Q_{s} -G_{ik}^{(r-1)}   )r_{ij} +  } \end{array}\]
\begin{equation} \label{eq56}
+      \sum _{q=1}^{r-1} B_{i} D_{ik}^{-1} B_{i}^{T} G_{ik}^{(r-q)}  ,
\end{equation}
\[r>1;      \tilde{A}_{ik} \equiv A_{i} -B_{i} D_{ik}^{-1} B_{i}^{T} G_{ik}^{(0)} ,      i=\overline{1,l},    k\geq0.\]

Note that the system (55) consists of independent matrix equations which, for fixed $i=1,2,...,l$, give a solution to the problem of optimal stabilization of the system
\begin{equation} \label{eq57}
dx(t)=(A_{i} x(t)+B_{i} u(t))dt+\sigma _{i} x(t)dw(t)+ C_i x(t)dN(t),
\end{equation}
with the quality criterion
\[{\it I}_{u} (y,h,x_0)=\sum _{k=0}^{\infty }\int _{t_{k} }^{\infty }{\it E}\{ x^{T} (t)M_{ik} x(t)+u^{T} (t)D_{ik} u(t)/ x_{0}\} dt ,\]
\begin{equation} \label{eq58}
i=\overline{1,l},   k\ge 0,     M_{ik} >0,  D_{ik} >0.
\end{equation}

A necessary and sufficient condition for the solvability of the system (55) is the existence of linear admissible control in the system (57), which provides exponential stability in the mean square of the unperturbed motion of this system \cite{LM10}.

Let us assume that the system of matrix quadratic equations (55) has a unique positive definite solution $G_{ik}^{(0)} >0,    i=\overline{1,l},    k\ge 0$.

Equation (56) to find $G_{ik}^{(r)} >0,  r\ge 1,  k\ge 0$ are linear, so they have a unique solution for fixed $i=\overline{1,l},    k\ge 0,    r\ge 1$ and any matrices that are on the right side of (56).

Indeed, the system
\begin{equation} \label{eq59}
dx(t)=\tilde{A}_{ik} x(t)dt+\sigma _{i} x(t)dw(t)+C_ix(t)dN(t)
\end{equation}
obtained by closure the system (57) with the optimal control
\[u_{k}^{0} =-D_{ik}^{-1} B_{ik}^{T} G_{ik}^{(0)} x(t),\]
which provides exponential stability in the mean square. Then there is a unique solution to the system (56). Note that in the linear case for autonomous systems, the asymptotic stability is equivalent to the exponential stability \cite{TYM11}. Consider a theorem, which originates from the results of this work.

\textbf{Theorem 4.} 
If a strong solution $x(t)$ of the system (57) is exponentially stable in the mean square, then there exist Lyapunov functions $v_{k} (y,h,x),  k\ge 0$, which satisfy the conditions:
\[c_{1} \left\| x\right\| ^{2} \le v_{k} (y,h,x)\le c_{1} \left\| x\right\| ;\]
\[\frac{d{\it E}[v_{k} ]}{dt} \le -c_{3} \left\| x\right\| ^{2} .\]

Thus, the system of matrix equations (55), (56) allows us to find consistently the coefficients $G_{ik}^{(r)} >0$ of the corresponding series (53), (54), starting with a positive solution $G_{ik}^{(0)} >0,  i=\overline{1,l}, k\ge 0$ of the system (55).

The next step is to prove the convergence of the series (53), (54).
Without the loss of generality, we simplify notations by fixing $k \geq 0$.
Denote $L_{r} := \mathop{\max }\limits_{\begin{array}{l} {i=\overline{1,l},} \\ {k\ge 0} \end{array}} \left\| G_{i}^{(r)} \right\| $. Then from (56) it follows that there is a constant $c>0$, that for any $r>0$ the following estimate is correct
\begin{equation} \label{eq60}
L_{r} \le c\left[\sum _{q=1}^{r-1}L_{q} L_{r-q}  +L_{r-1} \right].
\end{equation}

Next, we use the method of majorant series.

Consider the quadratic equation
\begin{equation} \label{eq61}
\rho ^{2} +(a+\varepsilon )\rho +b=0,
\end{equation}
where the coefficients $a,b$ are chosen such that the power series expansion of one of the roots of this equation is a majorant series for (53).

We have
\begin{equation} \label{eq62}
\rho _{1,2} =-\frac{a+\varepsilon }{2} \pm \sqrt{\frac{(a+\varepsilon )^{2} }{4} -b} =\sum _{r=0}^{\infty }\varepsilon ^{r} \rho _{r}  .
\end{equation}

Let us substitute (62) into (61), and equate coefficients at equal powers of $\varepsilon $. Then we get an expression for $\rho_r$ through $\rho _{0} ,...,\rho _{r-1} $:
\begin{equation} \label{eq63}
\rho _{r} =-\frac{1}{2\rho _{0} +a} \left[\sum _{q=1}^{r-1}\rho _{q} \rho _{r-q} +\rho _{r-1}  \right],
\end{equation}
where $\rho_0$ should be found from the equation
\begin{equation} \label{eq64}
\rho _{0}^{2} +a\rho _{0} +b=0.
\end{equation}

Comparing (60) and (63), we see that the series (62) will be major for (53), if we consider
\[c=-\frac{1}{2\rho _{0} +a} >0;    ~  \rho _{0} =L_{0} >0.\]

Thus, the values of the coefficients $a, b$ in the equation (61) are
\[\begin{array}{l} {a=-\left[\frac{1}{c} +2L_{0} \right]<0;      } \\ {b=\frac{L_{0} }{c} +L_{0}^{2} >0.} \end{array}\]

Using the known $a$ and $b$ from (61) we find that the majorant series for (53) will be the expansion one of the roots of (61). This root is such that its values are determined by
\[\rho _{0} =L_{0} =-\frac{a}{2} -\sqrt{\frac{a^{2} }{4} -b} .\]

Convergence of the series (53) for $v_{k}^{0} (y,h,x)$ follows from the obvious inequality
\[\left\| \sum _{r=0}^{\infty }\varepsilon ^{r} G^{(r)} (y,h) \right\| \le \sum _{r=0}^{\infty }L_{r} \varepsilon ^{r}  .\]

Thus, we have proved the assertion which is formalized below as the Theorem 5:

\textbf{Theorem 5.}  Let:
\begin{enumerate}
\item for $\forall i=\overline{1,l}, k\ge 0$ the system (57) has a linear admissible control;
\item transition intensities $q_{ij}$ of a homogeneous Markov chain $\xi$ satisfy the condition (51).
\end{enumerate}

Then
\begin{enumerate}
\item  there is a unique solution to the problem of optimal stabilization of the system (43), (39), (40) with the jump condition (41) of the phase vector $x\in {\it R}^{m} $;
\item  the optimal Lyapunov function $v_{k}^{0} (y,h,x)$ and optimal control $u_{k}^{0} (y,h,x)$ determined by convergent series (53), (54), whose coefficients are found from the corresponding systems (55), (56).
\end{enumerate}

\textbf{Case II.} Let us substitute the series $G_{ik} =\sum _{r=0}^{\infty }\varepsilon ^{r} G_{ik}^{(r)}$ in (44) and equate the coefficients at the same powers $\varepsilon$. Then, taking into account (52), we obtain the following equations:
\[G_{ik}^{(0)} A_{i} +A_{i}^{T} G_{ik}^{(0)} +\sigma _{i}^{T} G_{ik}^{(0)} \sigma _{i} - B_{i} D_{ik}^{-1} B^T_{i} G_{ik}^{(0)} +\]
\begin{equation} \label{eq65}
+\lambda C_i^T G_{ik}C_i+\sum _{j\ne i}^{l}(G_{jk}^{(0)} -G_{ik}^{(0)} )q_{ij} +M_{ik}  =0,  k\ge 0,
\end{equation}
\begin{equation} \label{eq66}
G_{ik}^{(r)} \tilde{A}_{ik} +\tilde{A}_{ik}^{T} G_{ik}^{(r)} +\sigma _{i}^{T} G_{ik}^{(r)} \sigma _{i} +\lambda C_i^T G_{ik}C_i+\sum _{j\ne i}^{l}(G_{jk}^{(r)} -G_{ik}^{(r)} )q_{ij}  ={  \Phi }_{ik}^{(r)}   ,
\end{equation}
where $i=\overline{1,l}    k\ge 0$, $  \tilde{A}_{ik} = A_{i} -B_{i} D_{ik}^{-1} B_{i}^{T} G_{ik}^{(0)} ,$
$$
{  \Phi }_{ik}^{(r)} =\sum _{q=1}^{r-1} B_{i} D_{ik}^{-1} B_{i}^T G_{ik}^{(r-q)}  -
$$
\[-\sum _{j\ne i}^{l}(K_{ij}^{T} G_{jk}^{(r-1)} +G_{jk}^{(r-1)} K_{ij} +K_{ij}^{T} G_{jk}^{(r-2)} K_{ij} +\sum _{s=1}^{N}Q_{s}^{T} G_{jk}^{(r-2)} Q_{s}   )q_{ij}. \]

Based on the above-mentioned, the following theorem is correct:

\textbf{Theorem 6.} Let

\begin{enumerate}
\item  the system of matrix equations (65) has a unique positive definite solution $G_{ik}^{(0)} >0,    i=\overline{1,l};    k\ge 0$;
\item jumps of the phase vector $x\in {\it R}^{m} $ satisfy the condition (52).
\end{enumerate}

Then the linear-quadratic optimal stabilization problem (43), (39), (40) of minimizing the functional (42) has a unique solution, which is given in the form of convergent series (53), (54), and the matrices  $G_{ik}^{(r)} ,    i=\overline{1,l};    r\ge 1,    k\ge 0$ is the only solution to the linear matrix equations (66).

\section{Discussion}

In this work, we have obtained sufficient conditions for the existence of an optimal solution for a stochastic dynamical system with jumps, which transform the system to a stable one in probability. The second Lyapunov method was used to investigate the existence of an optimal solution. This method is efficient both for ordinary differential equations (ODE) and for stochastic differential equations (SDE). As it can be seen from the proof of the Theorem 2, the existence of finite bounds for jumps at non-random time moments $t_m$ ($\lim_{m\to \infty }t_m =T^*<\infty$) does not impact the stability of the solution. On the other hand, the fact that $|t_{m+} - t_m|>\delta, m\ge 1$ was used for proving the existence of the optimal control (Theorem 3). This restriction is also present in the works of other authors. Thus, the goal of the subsequent work could be to construct an optimal control without the assumption $|t_{m+} - t_m|>\delta, m\ge 1$, which will considerably expand the scope of the second Lyapunov method.

\section{Conclusions}

In this work, we obtained sufficient conditions for the existence of a solution to an optimal stabilization problem for dynamical systems with jumps. We considered the case of a linear system with a quadratic quality functional. We showed that the problem of designing an optimal control that stabilizes the system to a stable one in probability, reduces to the problem of solving the Riccati equations. Additionally, for a linear autonomous system, the method of a small parameter is substantiated for solving the problem of optimal stabilization. The obtained solutions can be used to describe a model of a stock market in economics, biological systems including models of cancer response to treatment and other complex dynamical systems.\newline

\section*{Acknowledgments}
We would like to acknowledge the administrations of the Luxembourg Institute of Health (LIH) and Luxembourg National Research Fund (FNR) for their support in organizing scientific contacts between research groups in Luxembourg and Ukraine.

\end{document}